\documentclass[11pt]{article}

\usepackage[a4paper,margin=1in]{geometry}
\usepackage{amsmath,amssymb,amsthm,mathtools,bm}
\usepackage{enumitem}
\usepackage{hyperref}
\usepackage{mathrsfs}
\usepackage{array}
\usepackage{booktabs}
\usepackage{xcolor}

\usepackage{authblk} 
\usepackage[numbers,sort&compress]{natbib} 

\usepackage{tikz}
\usetikzlibrary{
  arrows.meta,
  calc,
  positioning,
  decorations.pathreplacing,
  patterns,
  fit,
  backgrounds,
  intersections,
  shapes.geometric  
}

\hypersetup{
  colorlinks=true,
  linkcolor=blue!60!black,
  citecolor=blue!60!black,
  urlcolor=blue!60!black
}

\newtheorem{theorem}{Theorem}[section]
\newtheorem{proposition}[theorem]{Proposition}

\theoremstyle{definition}
\newtheorem{definition}[theorem]{Definition}
\newtheorem{remark}[theorem]{Remark}

\newtheorem{assumption}[theorem]{Assumption}
\newtheorem{exercise}{Exercise}[section]
\newtheorem{solution*}{Solution}[section]
\newcommand{\R}{\mathbb{R}}

\newcommand{\Gr}{\operatorname{Gr}}

\newcommand{\argmin}{\operatorname*{arg\,min}}

\newcommand{\Tcone}{\mathcal{T}}
\newcommand{\VI}{\mathrm{VI}}
\newcommand{\Fset}{\mathcal{F}}
\newcommand{\Lcone}{\mathcal{L}}

\tikzset{
  axis/.style={->, thick},
  ray/.style={-Latex, thick},
  mybox/.style={draw, rounded corners, thick, inner sep=6pt, align=center},
  myellipse/.style={draw, ellipse, thick, align=center, inner sep=4pt},
  point/.style={circle, fill=black, inner sep=1.6pt},
  critical/.style={very thick, dashed},
  feasible/.style={thick},
  labelbox/.style={draw, rounded corners, fill=gray!8, inner sep=4pt, align=left}
}

\title{Optimization Workshop Notes for Mathematical Programming with Equilibrium Constraints (MPECs): \\ Verification of MPEC Hypotheses}
\author{Jiguang Yu\thanks{Email: jyu678@bu.edu}}
\affil{College of Engineering, Boston University, Boston, 02215, MA, USA}
\begin{document}

\maketitle

\begin{abstract}
In this workshop, we present a compact but rigorous introduction to the basic language
of nonlinear programming, variational inequalities, and complementarity systems.
The goal is twofold.
First, we explain the mathematical logic of hypotheses under which first-order optimality conditions for MPECs
become valid.
Second, we explain how to use that theory in research practice:
how to classify a model, choose the appropriate verification route, prove the right
hypotheses, and derive a correct first-order analysis.
\end{abstract}

\tableofcontents

\section{Problem Setting and Motivation}

We study the MPEC
\begin{equation}\label{eq:MPEC}
\begin{aligned}
\min_{x,y} \quad & f(x,y) \\
\text{s.t.} \quad & x \in X, \\
& y \in S(x),
\end{aligned}
\end{equation}
where
\begin{enumerate}[label=(\roman*)]
    \item $f:\R^{n+m}\to\R$ is continuously differentiable,
    \item $X\subset \R^n$ is a closed convex set,
    \item $S(x)$ is the solution set of a parametric variational inequality
    \[
    \VI(F(x,\cdot),C(x)).
    \]
\end{enumerate}
That is,
\[
S(x):=\left\{y\in C(x)\;:\; (v-y)^\top F(x,y)\ge 0 \quad \forall v\in C(x)\right\}.
\]

Define the feasible set
\[
\Gamma:=\{(x,y)\in X\times \R^m : y\in S(x)\}.
\]

If, near a reference point $(\bar x,\bar y)\in \Gamma$, the lower-level solution is
locally unique and can be written as a function
\[
y = y(x),
\]
then the MPEC \eqref{eq:MPEC} reduces locally to
\begin{equation}\label{eq:reduced}
\min_{x\in X} \ \varphi(x):=f(x,y(x)).
\end{equation}
This is attractive for three reasons:
\begin{enumerate}[label=(\roman*)]
    \item \textbf{Economic interpretation:} the follower's equilibrium response becomes a locally well-defined reaction map.
    \item \textbf{Geometric simplification:} the tangent cone to $\Gamma$ can be expressed in terms of the derivative of $y(\cdot)$.
    \item \textbf{Algorithmic relevance:} local descent methods can be designed on the reduced problem \eqref{eq:reduced}, provided one can evaluate or approximate $y'(x;d)$.
\end{enumerate}

The cost of this simplification is that one needs strong local stability properties of the lower-level variational inequality.

\subsection{The Lower-Level Variational Inequality}

Assume that, near $\bar x$, the feasible set of the lower-level problem has the form
\[
C(x)=\{y\in \R^m : g(x,y)\le 0\},
\]
where $g=(g_1,\dots,g_\ell):\R^{n+m}\to \R^\ell$ is $C^2$, and for each fixed $x$ the functions
$g_i(x,\cdot)$ are convex in $y$.

The lower-level variational inequality at parameter $x$ is:

\medskip
\noindent
\textbf{VI$(F(x,\cdot),C(x))$:}
find $y\in C(x)$ such that
\[
(v-y)^\top F(x,y)\ge 0 \qquad \forall v\in C(x).
\]

We will focus on a given feasible pair $(\bar x,\bar y)\in \Gamma$ satisfying
\[
\bar y \in S(\bar x).
\]

Under suitable constraint qualifications, the solution $\bar y$ admits Lagrange multipliers
$\lambda\in \R^\ell_+$ satisfying
\begin{equation}\label{eq:kkt-lower}
F(\bar x,\bar y)+\sum_{i=1}^{\ell}\lambda_i \nabla_y g_i(\bar x,\bar y)=0,
\end{equation}
together with complementarity
\[
\lambda_i\ge 0,\qquad g_i(\bar x,\bar y)\le 0,\qquad \lambda_i g_i(\bar x,\bar y)=0
\quad (i=1,\dots,\ell).
\]

Let the active set be
\[
I(\bar x):=\{i : g_i(\bar x,\bar y)=0\}.
\]

The multiplier set is
\[
\mathcal M(\bar x):=
\left\{
\lambda\in \R^\ell_+:
F(\bar x,\bar y)+\sum_{i=1}^{\ell}\lambda_i\nabla_y g_i(\bar x,\bar y)=0,\ 
\lambda_i g_i(\bar x,\bar y)=0 \ \forall i
\right\}.
\]

\subsection{B-Differentiability and the Key Implicit Assumption}

The reduced problem \eqref{eq:reduced} is only useful if the reaction map $y(\cdot)$ has enough regularity.

\begin{definition}[Directional and B-differentiability]
Let $G:W\subset \R^p\to \R^q$ be defined on an open set $W$ and let $z\in W$.
\begin{enumerate}[label=(\roman*)]
    \item $G$ is \emph{directionally differentiable} at $z$ if for every $d\in \R^p$,
    \[
    G'(z;d):=\lim_{t\downarrow 0}\frac{G(z+td)-G(z)}{t}
    \]
    exists.
    \item $G$ is \emph{B-differentiable} at $z$ if it is locally Lipschitz near $z$ and directionally differentiable at $z$.
\end{enumerate}
\end{definition}

\begin{remark}
B-differentiability is the natural regularity notion in this context. It is weaker than Fr\'echet differentiability but strong enough to support tangent-cone calculations and chain rules for directional derivatives.
\end{remark}

\begin{assumption}[Basic implicit function property]\label{ass:BIF}
At $(\bar x,\bar y)\in \Gamma$, there exist neighborhoods $U$ of $\bar x$ and $V$ of $\bar y$, and a Lipschitz continuous function
\[
y:U\cap X \to V
\]
such that
\begin{enumerate}[label=(\roman*), leftmargin=2em]
    \item $y(\bar x)=\bar y$,
    \item $y(\cdot)$ is directionally differentiable at $\bar x$,
    \item for each $x\in U\cap X$, the point $y(x)$ is the unique solution in $V$ of $\VI(F(x,\cdot),C(x))$.
\end{enumerate}
We call this assumption \textbf{BIF}.
\end{assumption}
Under BIF, the feasible set $\Gamma$ is locally the graph of the map $y(\cdot)$. Therefore its tangent cone at $(\bar x,\bar y)$ becomes
\[
\Tcone((\bar x,\bar y);\Gamma)
=
\left\{(dx,dy)\in \Tcone(\bar x;X)\times \R^m : dy = y'(\bar x;dx)\right\}.
\]

This identity is the central geometric gain of the implicit approach.

\subsection{Stationarity of the MPEC Under BIF}

Suppose BIF holds at $(\bar x,\bar y)$ and let
\[
\varphi(x):=f(x,y(x)).
\]

Then the chain rule for directional derivatives yields
\[
\varphi'(\bar x;dx)=\nabla_x f(\bar x,\bar y)^\top dx
+\nabla_y f(\bar x,\bar y)^\top y'(\bar x;dx).
\]

Hence $(\bar x,\bar y)$ is stationary for the original MPEC if and only if
\begin{equation}\label{eq:stationarity-implicit}
\nabla_x f(\bar x,\bar y)^\top dx
+\nabla_y f(\bar x,\bar y)^\top y'(\bar x;dx)\ge 0
\qquad
\forall dx\in \Tcone(\bar x;X).
\end{equation}

\begin{remark}
Equation \eqref{eq:stationarity-implicit} is the exact analogue of first-order stationarity for the reduced problem \eqref{eq:reduced}. The main unresolved issue is now:
\[
\textit{How do we compute } y'(\bar x;dx)\textit{?}
\]
The answer is: by linearizing the lower-level VI, provided BIF can be justified.
\end{remark}

\subsection{How to Verify BIF: The Main Structural Assumptions}

We now introduce the regularity conditions guaranteeing BIF.

\begin{assumption}[MFCQ at $(\bar x,\bar y)$]
There exists $v\in \R^m$ such that
\[
\nabla_y g_i(\bar x,\bar y)^\top v < 0
\qquad \forall i\in I(\bar x).
\]
\end{assumption}
This ensures metric regularity of the active inequality system and nonemptiness/boundedness of lower-level multipliers near the reference point.

\begin{assumption}[Constant rank CQ at $(\bar x,\bar y)$]
There exists a neighborhood $W$ of $(\bar x,\bar y)$ such that, for every subset
$J\subset I(\bar x)$, the rank of
\[
\{\nabla_y g_i(x,y): i\in J\}
\]
is constant for $(x,y)\in W$.
\end{assumption}
\begin{remark}
CRCQ is weaker than LICQ and allows nonunique multipliers. It is precisely this flexibility that makes the theory substantially broader than classical strong regularity theory for parametric NLPs under LICQ.
\end{remark}

We introduce the SCOC family.
For every multiplier $\lambda\in \mathcal M(\bar x)$, define its support
\[
\mathrm{supp}(\lambda):=\{i:\lambda_i\neq 0\}.
\]

Consider the family $\mathcal B(\bar x)$ of index sets $J\subset I(\bar x)$ such that
\begin{enumerate}[label=(\roman*)]
    \item there exists $\lambda\in \mathcal M(\bar x)$ with $\mathrm{supp}(\lambda)\subset J$,
    \item the vectors $\{\nabla_y g_i(\bar x,\bar y): i\in J\}$ are linearly independent.
\end{enumerate}

Write
\[
\mathcal B(\bar x)=\{J_1,\dots,J_N\}.
\]

For each $J_j\in \mathcal B(\bar x)$, let $\lambda^j$ be the associated extreme multiplier and define
\[
L(x,y,\lambda):=F(x,y)+\sum_{i=1}^{\ell}\lambda_i \nabla_y g_i(x,y),
\]
together with the matrix
\begin{equation}\label{eq:Aj}
A_j:=
\begin{bmatrix}
\nabla_y L(\bar x,\bar y,\lambda^j) & \nabla_y g_{J_j}(\bar x,\bar y)^\top \\
-\nabla_y g_{J_j}(\bar x,\bar y) & 0
\end{bmatrix}.
\end{equation}

\begin{assumption}[Strong coherent orientation condition]
The matrices $A_1,\dots,A_N$ all have the same nonzero determinantal sign:
\[
\mathrm{sgn}(\det A_1)=\cdots = \mathrm{sgn}(\det A_N)\in \{-1,+1\}.
\]
\end{assumption}

\begin{remark}
SCOC is the decisive nondegeneracy condition in the theory. Under LICQ and unique multipliers, it collapses to the nonsingularity of a single KKT matrix. In general, it allows multiple multipliers while enforcing a global consistency of local orientations.
\end{remark}

\subsection{Piecewise Smooth Analysis and Why It Enters}

The lower-level VI is handled through its \emph{normal map}. Let
\[
\bar v := \bar y - F(\bar x,\bar y).
\]

Define the Euclidean projector onto $C(x)$ by
\[
\Pi_{C(x)}(v):=\argmin_{y\in C(x)} \frac12 \|y-v\|^2,
\]
and define the normal map
\begin{equation}\label{eq:normal-map}
H(x,v):=F(x,\Pi_{C(x)}(v))+v-\Pi_{C(x)}(v).
\end{equation}

Then the standard equivalence is:
\[
y \text{ solves } \VI(F(x,\cdot),C(x))
\quad \Longleftrightarrow \quad
v:=y-F(x,y) \text{ solves } H(x,v)=0,
\]
with
\[
y=\Pi_{C(x)}(v).
\]

The projector $(x,v)\mapsto \Pi_{C(x)}(v)$ is generally not $C^1$.
However, under MFCQ and CRCQ it is locally piecewise $C^1$. Since $F$ is $C^1$, the normal map $H$ is also PC$^1$.
This is important because PC$^1$ functions are locally Lipschitz, directionally differentiable, and amenable to inverse/implicit function theorems based on coherent orientation.

\subsection{Main Existence Theorem for the Implicit Solution Map}

We are now ready for the central theorem.

\begin{theorem}[Existence of a PC$^1$ implicit solution map]
Assume:
\begin{enumerate}[label=(\roman*)]
    \item $g$ is $C^2$ and each $g_i(x,\cdot)$ is convex in $y$,
    \item $F$ is $C^1$,
    \item $\bar y$ solves $\VI(F(\bar x,\cdot),C(\bar x))$,
    \item MFCQ, CRCQ, and SCOC hold at $(\bar x,\bar y)$.
\end{enumerate}
Then there exist neighborhoods $U$ of $\bar x$ and $V$ of $\bar y$, and a PC$^1$ function
\[
y:U\to V
\]
such that
\begin{enumerate}[label=(\alph*), leftmargin=2em]
    \item $y(\bar x)=\bar y$,
    \item for every $x\in U$, the point $y(x)$ is the unique solution in $V$ of
    $\VI(F(x,\cdot),C(x))$.
\end{enumerate}
In particular, BIF holds at $(\bar x,\bar y)$.
\end{theorem}

The proof proceeds through four conceptual steps:

\begin{enumerate}[label=(\roman*)]
    \item Show that $(x,v)\mapsto \Pi_{C(x)}(v)$ is PC$^1$ near $(\bar x,\bar v)$.
    \item Deduce that the normal map $H$ in \eqref{eq:normal-map} is PC$^1$.
    \item Prove that the partial directional derivative
    \[
    H'_v(\bar x,\bar v;\cdot)
    \]
    is invertible; this is where SCOC enters.
    \item Apply a PC$^1$ implicit function theorem to obtain a locally unique PC$^1$ map
    \[
    v=v(x) \quad \text{solving } H(x,v)=0,
    \]
    and then define
    \[
    y(x):=\Pi_{C(x)}(v(x)).
    \]
\end{enumerate}

\subsection{Computing the Directional Derivative of the Response Map}

The existence theorem guarantees BIF. We now derive $y'(\bar x;dx)$.

\subsubsection{Critical cone associated with a multiplier}

Given $\lambda\in \mathcal M(\bar x)$ and a direction $dx\in \R^n$, define the directional critical set
\begin{equation}\label{eq:critical-cone}
\mathcal K(\bar x,\lambda;dx)
:=
\left\{
dy\in \R^m :
\begin{array}{l}
\nabla_x g_i(\bar x,\bar y)^\top dx + \nabla_y g_i(\bar x,\bar y)^\top dy \le 0
\quad \text{if } i\in I(\bar x),\ \lambda_i=0,\\[0.3em]
\nabla_x g_i(\bar x,\bar y)^\top dx + \nabla_y g_i(\bar x,\bar y)^\top dy = 0
\quad \text{if } \lambda_i>0
\end{array}
\right\}.
\end{equation}

Define also the critical multiplier set
\[
\mathcal M^c(\bar x;dx):=
\mathrm{arg\,max}_{\lambda\in \mathcal M(\bar x)}
\sum_{i=1}^{\ell}\lambda_i \nabla_x g_i(\bar x,\bar y)^\top dx.
\]

\subsubsection{The linearized variational inequality}

For a multiplier $\lambda\in \mathcal M^c(\bar x;dx)$, define the affine variational inequality
\begin{equation}\label{eq:AVI}
\text{find } dy\in \mathcal K(\bar x,\lambda;dx)
\text{ such that }
(w-dy)^\top
\Big(
\nabla_x L(\bar x,\bar y,\lambda)dx
+
\nabla_y L(\bar x,\bar y,\lambda)dy
\Big)
\ge 0
\quad \forall w\in \mathcal K(\bar x,\lambda;dx).
\end{equation}

\begin{theorem}[Directional derivative formula]
Under the hypotheses of the main existence theorem, for every $dx\in \R^n$ and every
$\lambda\in \mathcal M^c(\bar x;dx)$, the directional derivative $y'(\bar x;dx)$ solves the AVI
\eqref{eq:AVI}. If $\lambda$ is an extreme multiplier, then the AVI has a unique solution.
\end{theorem}

\begin{remark}
Thus the derivative of the implicit response map is not obtained from a single linear system in general. It is obtained from an auxiliary variational inequality over a direction-dependent critical cone. This is the correct linearization when the lower-level problem has inequality constraints and potentially multiple multipliers.
\end{remark}

\subsection{Tangent Cone Representation and MPEC Constraint Qualifications}

Once $y'(\bar x;dx)$ is available, the tangent cone to the MPEC feasible set is fully characterized:
\[
\Tcone((\bar x,\bar y);\Gamma)
=
\left\{
(dx,dy)\in \Tcone(\bar x;X)\times \R^m :
dy = y'(\bar x;dx)
\right\}.
\]

This representation can then be combined with the directional derivative formula above to verify MPEC-specific constraint qualifications, including the extreme CQ in the framework of the text.

\begin{remark}
The notable point is that CRCQ is not merely technical decoration. It is essential for obtaining the piecewise smooth structure that allows one to prove the relevant tangent-cone identity and the ensuing first-order conditions.
\end{remark}

\subsection{When Is the Implicit Solution Map Actually Fr\'echet Differentiable?}

B-differentiability of $y(\cdot)$ is guaranteed by the theorem, but Fr\'echet differentiability is subtler.

Define
\[
\mathscr{L}(\bar x,dx)
:=
\left\{
dy\in \R^m :
\nabla_x g_i(\bar x,\bar y)^\top dx+\nabla_y g_i(\bar x,\bar y)^\top dy = 0
\quad \forall i\in I(\bar x)
\right\}.
\]

Then one has the following qualitative criterion:

\begin{theorem}[Characterization of Fr\'echet differentiability]
Under the hypotheses of the main existence theorem, the following are equivalent:
\begin{enumerate}[label=(\roman*), leftmargin=2em]
    \item $y(\cdot)$ is Fr\'echet differentiable at $\bar x$,
    \item for every $dx\in \R^n$, the directional derivative $y'(\bar x;dx)$ lies in $\mathscr{L}(\bar x,dx)$,
    \item equivalently, the AVI defining $y'(\bar x;dx)$ reduces to a linear system on the affine subspace of equality-consistent directions.
\end{enumerate}
\end{theorem}

\begin{remark}
This theorem makes a key conceptual point: directional differentiability is generic in the PC$^1$ setting, but full Fr\'echet differentiability requires that the one-sided directional laws glue together linearly across all directions.
\end{remark}

\subsection{The LICQ Case: Cleaner Structure}

If the lower-level problem satisfies LICQ at $(\bar x,\bar y)$, then the multiplier is unique:
\[
\mathcal M(\bar x)=\{\bar\lambda\}.
\]

In this case,
\begin{enumerate}[label=(\roman*)]
    \item the implicit solution map $y(\cdot)$ exists and is PC$^1$ under the theorem's assumptions,
    \item there is also a PC$^1$ multiplier map $\lambda(x)$,
    \item the derivative pair $(y'(\bar x;dx),\lambda'(\bar x;dx))$ is obtained from a linearized KKT system.
\end{enumerate}
This is the closest analogue to the classical sensitivity theory of nonlinear programming, but the current framework is broader because it does not require LICQ in the general theorem.

\subsection{Interpretation of SCOC}

SCOC is best understood as a generalized strong regularity condition.

If LICQ holds and the multiplier is unique, SCOC reduces to the nonsingularity of the single KKT matrix
\[
\begin{bmatrix}
\nabla_y L(\bar x,\bar y,\bar\lambda) & \nabla_y g_{I(\bar x)}(\bar x,\bar y)^\top\\
-\nabla_y g_{I(\bar x)}(\bar x,\bar y) & 0
\end{bmatrix}.
\]

When multipliers are not unique, different active-index selections can produce different local linear models. SCOC requires that all such models have the same nonzero orientation. This coherence prevents local folds and branching pathologies in the implicit equation for the normal map.

\begin{remark}
So the phrase \emph{strong coherent orientation} is not cosmetic: it literally combines
\begin{enumerate}[label=(\roman*)]
    \item a \emph{strong} nondegeneracy requirement (nonzero determinant),
    \item and a \emph{coherence} requirement (same sign across all admissible pieces).
\end{enumerate}
\end{remark}

Despite its elegance, the IMP approach is restrictive.
\begin{enumerate}[label=(\roman*)]
    \item It requires a local single-valued solution branch $y(x)$.
    \item The feasible region must have the structural form $X\times \R^m$ at the upper level.
    \item The regularity assumptions can be unnecessarily strong for some MPECs.
\end{enumerate}
This motivates alternative piecewise or set-valued approaches in which one works directly with multiple local response branches or with more flexible linearizations of the feasible set.

\subsection{Takeaways}

For research purposes, the section can be summarized as follows.
\begin{enumerate}[label=(\roman*)]
    \item The implicit programming approach turns an MPEC into a reduced nonsmooth optimization problem whenever a locally unique response map exists.
    \item The response map is typically not globally smooth; the correct local regularity class is PC$^1$.
    \item MFCQ + CRCQ + SCOC are the structural assumptions that deliver local uniqueness and PC$^1$ regularity.
    \item The derivative of the response map is computed by a direction-dependent affine variational inequality.
    \item These results provide the tangent-cone formula and hence the first-order stationarity conditions for the MPEC.
\end{enumerate}

The implicit programming approach is a sensitivity-theoretic route to MPEC analysis.
Its real achievement is not merely to replace $y$ by $y(x)$, but to prove that such a replacement is mathematically legitimate in a highly nonsmooth setting involving variational inequalities, multiple multipliers, and nonpolyhedral feasible sets.

The main conceptual lesson is this:

\centerline{\emph{MPEC analysis is fundamentally an exercise in stability of generalized equations.}}

The implicit function $y(x)$ exists only when the lower-level equilibrium is sufficiently stable, and the language best suited to this stability is the language of
\begin{enumerate}[label=(\roman*)]
    \item piecewise smooth maps,
    \item coherent orientation,
    \item tangent cones,
    \item and affine variational inequalities.
\end{enumerate}

Those tools are what make the reduced formulation rigorous rather than heuristic.

\section{Part A. Theory from Chapter 4}

A mathematical program with equilibrium constraints (MPEC) is an optimization problem
in which part of the feasible set is defined by an equilibrium condition, typically a
variational inequality, complementarity system, or lower-level KKT system.
This structure appears naturally in engineering design, equilibrium modeling, and
multilevel games.
It is mathematically difficult because the feasible set is often nonsmooth, nonconvex,
and may even be disconnected. 

A generic form is
\[
\min_{(x,y)\in Z} f(x,y)
\qquad \text{subject to} \qquad
y \in S(x),
\]
where \(x\) is the upper-level variable, \(y\) is the lower-level response,
\(Z \subseteq \R^{n+m}\) is an explicit constraint set, and \(S(x)\) is the solution set
of a parametric variational inequality or complementarity problem.

The central issue is not merely how to rewrite the equilibrium constraint.
The central issue is this:

\begin{quote}
How do we verify assumptions under which first-order optimality conditions are actually
valid for the resulting MPEC?
\end{quote}

That is the purpose of Chapter 4, and that is the focus of these notes.

\subsection{Notation and basic viewpoint}

We work with feasible points \(z=(x,y)\in \R^{n+m}\).
The feasible region of the MPEC is denoted by \(\Fset\).
At a feasible point \(z\), the tangent cone is denoted by \(\Tcone(z;\Fset)\),
and a suitable linearized cone by \(\Lcone(z;\Fset)\).

The philosophy of Chapter 4 is geometric:

\begin{quote}
If the correct linearized cone captures the true tangent cone,
then stationarity can be characterized by multiplier systems analogous
to KKT conditions.
\end{quote}

For MPECs, this is a nontrivial statement.
Ordinary NLP constraint qualifications cannot be assumed mechanically,
because complementarity relations and equilibrium constraints create
degeneracy and nonsmoothness.

In classical smooth nonlinear programming, one often proves KKT conditions from
standard constraint qualifications such as LICQ or MFCQ.
In an MPEC, that route is unreliable.
Complementarity constraints typically violate classical regularity assumptions
at any point where both sides of the complementarity system are active.

Therefore, one must use the internal structure of the equilibrium system itself.
Chapter 4 does exactly that.
It identifies cases where the local geometry is regular enough that first-order
conditions can still be justified rigorously.

\subsection{AVI-constrained MPECs}

The cleanest case is the affine variational inequality (AVI) constrained mathematical
program, also called the MPAEC.
The model has the form
\begin{equation}
\label{eq:mpaec}
\begin{aligned}
\min_{x,y} \quad & f(x,y) \\
\text{s.t.} \quad & (x,y)\in Z,\\
& Dx + Ey + b \le 0,\\
& (y'-y)^\top(Px+Qy+q) \ge 0
\quad \forall y' \text{ such that } Dx+Ey'+b \le 0.
\end{aligned}
\end{equation}

This problem is special because the lower-level equilibrium mapping is affine in
\((x,y)\), and the lower-level feasible set is polyhedral.

\subsubsection*{Equivalent KKT representation}

The inner AVI can be expressed through a multiplier vector \(\lambda \ge 0\):
\begin{equation}
\label{eq:avi-kkt}
\begin{aligned}
\min_{x,y,\lambda}\quad & f(x,y)\\
\text{s.t.}\quad & (x,y,\lambda)\in Z\times \R^\ell_+,\\
& Dx+Ey+b \le 0,\\
& Px+Qy+q+E^\top \lambda = 0,\\
& \lambda^\top(Dx+Ey+b)=0.
\end{aligned}
\end{equation}

Although this resembles a linearly constrained NLP, it is not one:
the complementarity relation is quadratic and induces a piecewise feasible geometry.

\subsubsection*{The decisive geometric result}

The key point of Section 4.1 is that, in the AVI case, the relevant tangent and
linearized cones coincide.
Conceptually, the first-order model is exact.

\begin{proposition}[Geometric exactness in the AVI case]
For an AVI-constrained MPEC with polyhedral structure, the tangent cone and the
appropriate linearized cone coincide:
\[
\Tcone(z;\Fset)=\Lcone(z;\Fset).
\]
Consequently, the relevant MPEC constraint qualifications hold.
\end{proposition}

This is one of the strongest and most useful conclusions in the chapter.
It means that stationarity conditions can be derived without adding an external
constraint qualification beyond the structural assumptions already built into the AVI model.

\subsubsection*{Why this is powerful}

The result depends on two features:

\begin{enumerate}[label=(\roman*)]
    \item the lower-level VI is affine, so the directional critical structure can be analyzed explicitly;
    \item the feasible set is piecewise polyhedral, so tangent directions can be realized by feasible arcs.
\end{enumerate}

This is why the AVI case behaves much better than the general nonlinear VI case.

\subsubsection*{Piecewise polyhedrality and local minimality}

When \(Z\) is polyhedral, the feasible set is locally a finite union of convex polyhedra.
That allows one to strengthen stationarity to local minimality under pseudoconvexity.

\begin{proposition}[Sufficiency under pseudoconvexity]
If \(Z\) is polyhedral and the objective \(f\) is pseudoconvex on the feasible region,
then every stationary point of the AVI-constrained MPEC is a local minimum.
\end{proposition}

This conclusion is local, not global.
Even in the AVI case, a stationary point need not be a global minimizer.

\subsection{The implicit programming approach}

The next route is the implicit programming (IMP) approach.
Here the feasible set has the special product form
\[
(x,y)\in (X\times \R^m)\cap \Gr(S),
\]
where \(X\subseteq \R^n\) is a manageable upper-level set and \(S(x)\) is the solution set
of a parametric variational inequality.

The MPEC is
\begin{equation}
\label{eq:imp-mpec}
\min_{(x,y)} f(x,y)
\qquad \text{s.t.} \qquad
(x,y)\in (X\times \R^m)\cap \Gr(S).
\end{equation}

\subsubsection*{Core idea}

If the lower-level equilibrium has a locally unique solution \(y(x)\) near a reference point,
then the MPEC can be locally reduced to
\begin{equation}
\label{eq:reduced}
\min_{x\in X} \tilde f(x):=f(x,y(x)).
\end{equation}

This is the conceptual appeal of IMP:
the follower's variable becomes a local response function of the leader's variable.

\subsubsection*{Why ordinary differentiability is too strong}

For parametric VIs, the solution map \(x\mapsto y(x)\) is often not classically smooth.
The correct notion is directional differentiability, usually together with local Lipschitz continuity.

\begin{definition}
Let \(G:W\subseteq \R^n\to \R^m\), and let \(\bar x\in W\).

\begin{enumerate}[label=(\roman*)]
    \item \(G\) is \emph{directionally differentiable} at \(\bar x\) if
    \[
    G'(\bar x;d):=\lim_{\tau\downarrow 0}\frac{G(\bar x+\tau d)-G(\bar x)}{\tau}
    \]
    exists for every \(d\in \R^n\).
    \item \(G\) is \emph{\(B\)-differentiable} at \(\bar x\) if it is locally Lipschitz near \(\bar x\)
    and directionally differentiable at \(\bar x\).
\end{enumerate}
\end{definition}

In Chapter 4, this level of smoothness is central because it is strong enough for first-order calculus,
but weak enough to be realistic for VI solution maps.

\subsubsection*{Interpretation}

The IMP route is therefore a sensitivity theory.
It asks whether the lower-level solution map is locally single-valued and directionally smooth enough
to justify first-order analysis of the reduced objective \(\tilde f(x)\).

\subsubsection*{Strengths and limitations}

The strengths are obvious:
\begin{enumerate}[label=(\roman*)]
    \item reduced-space analysis;
    \item a natural interpretation of \(y\) as the follower's response;
    \item possible algorithmic use in descent methods.
\end{enumerate}

The limitations are equally important:
\begin{enumerate}[label=(\roman*)]
    \item local single-valuedness may fail;
    \item the upper-level feasible set must have the product form \(X\times \R^m\);
    \item the required regularity can be stronger than necessary for some models.
\end{enumerate}

\subsection{The piecewise programming approach}

Section 4.3 develops a complementary viewpoint.
Instead of collapsing the lower-level response into a single implicit function,
one works directly with the branchwise structure of the feasible region.

This is often the right choice when active sets, complementarity regimes,
or multiplier patterns generate several local branches.
The feasible set is then understood piece by piece.

\begin{quote}
IMP is appropriate when the follower's response is locally single-valued and directionally smooth.
PCP is appropriate when the local geometry is intrinsically piecewise.
\end{quote}

The two methods are not interchangeable.
They are tailored to different structural regimes.

\subsection{Exact penalization}

The chapter ends with a first-order exact penalization result.
Conceptually, exact penalization means that, under suitable local assumptions,
the equilibrium constraint can be moved into the objective with a sufficiently large penalty parameter
without changing the local solutions.

This matters for two reasons:
\begin{enumerate}[label=(\roman*)]
    \item it provides another route to local optimality conditions;
    \item it motivates penalty-based numerical methods.
\end{enumerate}

Exactness is not automatic.
It rests on the same local regularity analysis that supports the IMP or PCP route.

\subsection{Main takeaways from Part A}

A senior research should retain the following ideas.

\begin{enumerate}[label=(\roman*)]
    \item MPEC first-order theory is fundamentally about local geometry and sensitivity analysis.
    \item The AVI case is exceptionally favorable because tangent and linearized cones coincide.
    \item The IMP route is a reduction principle based on local single-valuedness of the lower-level response.
    \item \(B\)-differentiability is often the correct first-order smoothness notion.
    \item No single method dominates: some problems are best treated implicitly, others piecewise.
\end{enumerate}

\section{Part B. How to Use the Theory in Research Practice}

\subsection{The practical question}

When you face a new MPEC, do not begin by writing multipliers.
Begin by answering a structural question:

\begin{quote}
Which local model from Chapter 4 actually matches my problem?
\end{quote}

That is the core research habit this chapter teaches.

\subsection{A workable research pipeline}

A reliable workflow is the following.

\begin{enumerate}[label=(\roman*)]
    \item \textbf{Write the lower-level system explicitly.}
    Put it in VI, complementarity, or KKT form.

    \item \textbf{Separate explicit and equilibrium constraints.}
    Decide what belongs to \(Z\) and what belongs to the solution map \(S(x)\).

    \item \textbf{Classify the lower-level structure.}
    Is it affine, nonlinear but locally regular, or genuinely piecewise?

    \item \textbf{Choose the Chapter 4 route.}
    AVI, IMP, PCP, or exact penalization.

    \item \textbf{Verify the hypotheses for that route.}
    This is the non-negotiable step.

    \item \textbf{Only then derive stationarity.}
    Do not write a KKT-type system before proving that the underlying feasible-direction model is valid.

    \item \textbf{State your conclusion precisely.}
    Distinguish among stationarity, local minimality, reduced necessary conditions,
    and exact penalty equivalence.
\end{enumerate}

\subsection{How to choose the route}

\subsubsection*{Route 1: AVI}

Use Section 4.1 when the lower-level system is affine and polyhedral:
\[
F(x,y)=Px+Qy+q, \qquad C(x)=\{y: Dx+Ey+b\le 0\}.
\]

This is the best-case regime.
You should aim to prove that tangent and linearized cones coincide and then derive stationarity directly.

\subsubsection*{Route 2: IMP}

Use Section 4.2 when you have good reason to believe that the lower-level equilibrium
defines a locally unique response \(y=y(x)\) near a point of interest.

Then your goal is to verify local single-valuedness and directional smoothness of \(y(x)\),
form the reduced objective \(\tilde f(x)=f(x,y(x))\), and work with the reduced problem.

\subsubsection*{Route 3: PCP}

Use Section 4.3 when the feasible region is naturally branchwise:
different active sets or complementarity regimes produce different local formulas,
and forcing everything into one implicit function would hide essential structure.

\subsection{The AVI checklist}

If your problem is AVI-constrained, the checklist is:

\begin{enumerate}[label=(\roman*)]
    \item Verify the lower-level map is affine.
    \item Verify the lower-level feasible set is polyhedral.
    \item Check whether the explicit upper-level set \(Z\) is polyhedral or at least tractable.
    \item Write the lower-level system in KKT/complementarity form.
    \item Invoke the tangent-equals-linearized result.
    \item Derive the stationarity system.
    \item If appropriate, use pseudoconvexity to strengthen stationarity to local minimality.
\end{enumerate}

\subsubsection*{What this gives you in a paper}

This route supports a strong statement of the form:

\begin{quote}
Because the lower-level equilibrium is affine and the local feasible geometry is piecewise polyhedral,
the tangent cone coincides with the linearized cone.
Hence the relevant first-order stationarity conditions hold without an additional CQ.
\end{quote}

That is a rigorous theorem, not a formal manipulation.

\subsection{The IMP checklist}

If you use the implicit-programming route, follow this order:

\begin{enumerate}[label=(\roman*)]
    \item Write the MPEC as
    \[
    (x,y)\in (X\times \R^m)\cap \Gr(S).
    \]
    \item Fix a reference point \((\bar x,\bar y)\).
    \item Prove local uniqueness of the lower-level solution:
    \[
    y=y(x)\quad \text{for } x \text{ near } \bar x.
    \]
    \item Verify the directional smoothness required by Chapter 4.
    \item Form the reduced objective
    \[
    \tilde f(x)=f(x,y(x)).
    \]
    \item Compute the directional derivative of \(\tilde f\).
    \item Apply first-order optimality conditions on \(X\).
    \item Translate the conclusion back to the original MPEC.
\end{enumerate}

\subsubsection*{Three common warnings}

\begin{enumerate}[label=(\roman*)]
    \item Never assume local single-valuedness just because you want a reduced model.
    \item Never use classical derivatives unless you have classical differentiability.
    \item Never claim global equivalence from a purely local sensitivity theorem.
\end{enumerate}

\subsection{The PCP checklist}

If the problem is intrinsically piecewise, proceed as follows.

\begin{enumerate}[label=(\roman*)]
    \item Identify local branches, active-set regimes, or complementarity patterns.
    \item Decompose the feasible region accordingly.
    \item Compute tangent directions on each relevant piece.
    \item Determine which pieces are actually reachable from the reference point.
    \item Assemble the total feasible-direction description.
    \item Derive stationarity by testing the directional derivative of the objective on all such directions.
\end{enumerate}

This route is often more robust than IMP when the response is multivalued.

\subsection{How to verify hypotheses efficiently}

A good habit is to verify assumptions in the following order:

\begin{enumerate}[label=(\roman*)]
    \item Geometry first.
    Is the feasible set polyhedral, piecewise polyhedral, or finitely branched?

    \item Smoothness second.
    Is the lower-level map affine, \(C^1\), piecewise \(C^1\), directionally differentiable,
    or \(B\)-differentiable?

    \item Sensitivity third.
    Is the solution map locally single-valued?
    Can its directional derivative be computed?

    \item Multipliers last.
    Are lower-level multipliers unique, nonunique, or branch-dependent?
\end{enumerate}

This order prevents a frequent mistake:
writing down a multiplier system before proving the feasible-direction model that justifies it.

\subsection{How theorem statements should be written}

A clean theorem in the style of Chapter 4 should look like this:

\begin{quote}
Assume the MPEC has structure A.
Assume the lower-level equilibrium satisfies properties B and C.
Then the local feasible geometry has property D.
Consequently, every local minimizer satisfies stationarity system E.
\end{quote}

Avoid vague statements such as:
\begin{quote}
``By standard MPEC theory, KKT-type conditions hold.''
\end{quote}
That is rarely precise enough.

\subsection{How to use the Theory in algorithmic work}

Although Chapter 4 is theoretical, it points directly to computational design.

If AVI applies: use the complementarity representation or branchwise polyhedral structure to justify
NLP reformulations, active-set methods, or SQP-style schemes.

If IMP applies: use the reduced function \(x\mapsto f(x,y(x))\) to derive reduced-space search directions.

If exact penalization applies: use a penalized model as a locally equivalent surrogate, while being explicit about the
local nature of that equivalence.

\subsection{The main mistakes to avoid}

\begin{enumerate}[label=\arabic*.]
    \item Do not confuse a reformulation with a proof.
    Writing the lower-level KKT conditions is not the same as verifying first-order theory.
    \item Do not import NLP constraint qualifications without checking whether complementarity destroys them.
    \item Do not assume more smoothness than your model actually has.
    \item Do not turn local sensitivity results into global claims.
    \item Do not ignore reachable-direction geometry.
\end{enumerate}

The right way to use Chapter 4 is not to memorize technical lemmas.
It is to internalize its logic:

\begin{quote}
Classify the lower-level structure.\\
Choose the correct local model.\\
Verify the hypotheses for that model.\\
Only then derive stationarity, sufficiency, or penalization.
\end{quote}

That logic is what turns an MPEC analysis into a rigorous mathematical argument.

\section{Test problems}

\subsection{True or False}

\begin{enumerate}[label=(\roman*)]
    \item For an Affine Variational Inequality (AVI) constrained mathematical program (MPAEC) with a polyhedral upper-level feasible set, the linearized cone $\mathcal{L}(\overline{z}; \mathcal{F})$ is strictly contained in the tangent cone $\mathcal{T}(\overline{z}; \mathcal{F})$ unless an additional Constraint Qualification (CQ), such as the Constant Rank Constraint Qualification (CRCQ), is explicitly satisfied.

\item For an MPAEC with a polyhedral set $Z$ and a pseudoconvex objective function $f$ on the feasible region $\mathcal{F}$, any stationary point $\overline{z}$ is guaranteed to be a global minimum.

\item If a function $G: W \subseteq \mathbb{R}^{n} \rightarrow \mathbb{R}^{m}$ is B-differentiable at $\overline{x}$, then $G$ is strongly F-differentiable at $\overline{x}$ if and only if its directional derivative $G^{\prime}(\overline{x}; d)$ is a linear function of the direction $d$.

\item If a $PC^{1}$ function $G: \mathbb{R}^{m} \rightarrow \mathbb{R}^{m}$ is coherently oriented at $\overline{v}$, its directional derivative $G^{\prime}(\overline{v}; \cdot)$ is automatically guaranteed to be a global Lipschitzian homeomorphism on $\mathbb{R}^{n}$.

\item Under the MFCQ, CRCQ, and SCOC assumptions for the lower-level variational inequality $(F(x, \cdot), C(x))$, the parametric projection mapping $\Pi_{C(x)}(v)$ is a $PC^{1}$ function near $(\overline{x}, \overline{v})$, and its $C^{1}$ pieces are determined by the equality-constrained nonlinear programs corresponding to the SCOC family of active index sets $\mathcal{B}(\overline{x})$.

\item In the specific case where the pair $(\overline{x}, \overline{y})$ is strongly nondegenerate for the lower-level variational inequality, the Strong Coherent Orientation Condition (SCOC) strictly requires the matrix $U^{T} \nabla_{y}L(\overline{x}, \overline{y}, \overline{\lambda}) U$ to be positive definite, where the columns of $U$ form an orthonormal basis for the null space of the active constraint gradients.

\item Assuming the conditions of Theorem 4.2.16 hold, the implicit solution function $y(x)$ of the parametric variational inequality is F-differentiable at $\overline{x}$ if and only if for all directions $dx \in \mathbb{R}^{n}$, the directional derivative $y^{\prime}(\overline{x}; dx)$ belongs to the affine subspace $\mathcal{E}(\overline{x}, dx)$defined by the gradients of all active constraints at $(\overline{x}, \overline{y})$.

\item In the Implicit Programming (IMP) approach, the Constant Rank Constraint Qualification (CRCQ) is an absolute necessity to establish the basic constraint qualification; without the CRCQ, it is fundamentally impossible to relate the tangent cone $\mathcal{T}(\overline{x}; \mathcal{F})$ to any linearized cone $\mathcal{L}^{\prime}(\overline{x}; \mathcal{F})$.

\item In the Piecewise Programming (PCP) approach, if the Strict Mangasarian-Fromovitz Constraint Qualification (SMFCQ) holds for the relaxed nonlinear program (NLP) at a local minimizer, it is mathematically guaranteed that the standard MFCQ will also hold for all corresponding local MPEC pieces.

\item Assuming the hypotheses for the existence of the $PC^{1}$ implicit function hold, and further assuming the upper-level set $X$ is compact and the lower-level feasible set $C(x)$ is bounded for each $x \in X$, there exists a finite penalty parameter $\overline{\rho} > 0$ such that the global minimizers of the normal-form MPEC are identical to the global minimizers of the exact penalty problem of order 1 for any $\rho \ge \overline{\rho}$.
\end{enumerate}

\subsubsection{Solutions for True or False}

\begin{enumerate}[label=(\roman*)]
    \item False. According to Proposition 4.1.1, for an MPAEC, $\mathcal{T}(\overline{z}; \mathcal{F}) = \mathcal{L}(\overline{z}; \mathcal{F})$ holds intrinsically. The text emphasizes that the extreme CQ must hold "without any particular assumptions," making the first-order KKT conditions necessary for local optimality without requiring additional CQs.
    \item False. While Proposition 4.1.4 proves that such a stationary point must be a local minimum, Example 4.1.6 explicitly demonstrates an AVI constrained MP with a convex objective function where a stationary point is a local minimum but not a global minimum.
    \item False. Based on Proposition 4.2.2(d), $G^{\prime}(\overline{x}; d)$ being linear in $d$ is the necessary and sufficient condition for $G$ to be F-differentiable at $\overline{x}$, not strongly F-differentiable. Strong F-differentiability requires a uniform limit condition for the directional derivatives in a neighborhood of $\overline{x}$(Proposition 4.2.2(e)).
    \item False. The text immediately following Lemma 4.2.13 states: "An example given in [85] shows that a $PC^{1}$function $G$ can be coherently oriented at a point $\overline{v}$ and yet none of the conditions (a), (b), or (c) in Lemma 4.2.13 holds." Coherent orientation is a prerequisite for these conditions, but it does not imply them without the additional assumption that the directional derivative is invertible.
    \item True. This is precisely the finding of Lemma 4.2.17. The proof demonstrates that extreme points of the KKT system dictate that the mapping relies on active index sets belonging to $\mathcal{B}(\overline{x})$, thereby defining the $C^{1}$ pieces of the $PC^{1}$ projection map.
    \item False. The text states that for strongly nondegenerate pairs, SCOC is equivalent to the nonsingularity of the matrix $U^{T} \nabla_{y}L(\overline{x}, \overline{y}, \overline{\lambda}) U$. It does not require positive definiteness, only that the determinant is non-zero.
    \item True. Theorem 4.2.28 explicitly establishes the equivalence between statement (a) (the implicit solution function $y(x)$ is F-differentiable at $\overline{x}$) and statement (b) (for all $dx \in \mathbb{R}^{n}$, $y^{\prime}(\overline{x}; dx) \in \mathcal{E}(\overline{x}, dx)$).
    \item False. The text notes in Section 4.2.6 that while the CRCQ guarantees the extreme CQ holds, "Without the CRCQ, although it might still be possible for the basic CQ to hold for the MPEC, the set of multipliers $M^{\prime}$ needed for equating the tangent cone... becomes difficult to identify." The text explicitly references Example 3.4.2 as a case where this set was successfully identified without CRCQ, meaning it is not "fundamentally impossible," just difficult.
    \item False. The text explicitly refutes this in Section 4.3.1: "The first example (which has appeared previously) illustrates that the SMFCQ holding for (11) [the relaxed NLP] does not necessarily imply that the MFCQ will hold for (3) [the local MPEC pieces]," and uses the revisited Example 4.3.4 to prove this failure.
    \item True. This is a direct summary of Theorem 4.4.3. The compactness of the feasible region and the Lipschitzian homeomorphism of the parametric normal map guarantee that the exact penalty function of order 1 shares the exact same global minimizers for a sufficiently large, finite penalty parameter.
\end{enumerate}

\subsection{Warm-up: geometric thinking for MPECs}

\begin{exercise}[Why standard NLP intuition fails]
Consider the complementarity set
\[
\Omega := \{(u,v)\in\R^2 : u\ge 0,\ v\ge 0,\ uv=0\}.
\]
\begin{enumerate}[label=(\alph*)]
    \item Sketch the set \(\Omega\).
    \item Compute the tangent cone \(\Tcone((0,0);\Omega)\).
    \item Show that \(\Omega\) is not a smooth manifold near the origin.
    \item Explain why any attempt to apply a standard NLP constraint qualification at \((0,0)\) is problematic.
\end{enumerate}
\end{exercise}

\begin{solution*}
\textbf{(a)} The condition \(uv=0\) together with \(u\ge 0\), \(v\ge 0\) implies that
either \(u\ge 0\) and \(v=0\), or \(u=0\) and \(v\ge 0\).
Hence
\[
\Omega = \{(u,0):u\ge 0\}\cup \{(0,v):v\ge 0\}.
\]
It is the union of the two nonnegative coordinate rays.

\textbf{(b)} By definition, \(d=(d_u,d_v)\in \Tcone((0,0);\Omega)\) if there exist
\((u_k,v_k)\in \Omega\) and \(t_k\downarrow 0\) such that
\[
\frac{(u_k,v_k)-(0,0)}{t_k}\to (d_u,d_v).
\]
Take first points on the horizontal branch \((u_k,0)\) with \(u_k\ge 0\):
this yields all directions \((\alpha,0)\) with \(\alpha\ge 0\).
Taking points on the vertical branch \((0,v_k)\) with \(v_k\ge 0\) yields all directions
\((0,\beta)\) with \(\beta\ge 0\).
No direction with both components strictly positive can occur, because every feasible point has
at least one zero component.
Thus
\[
\Tcone((0,0);\Omega)=\{(\alpha,0):\alpha\ge 0\}\cup \{(0,\beta):\beta\ge 0\}.
\]

\textbf{(c)} A smooth one-dimensional manifold near a point must have a unique tangent line there.
At the origin, \(\Omega\) has two distinct tangent rays, one horizontal and one vertical,
and no neighborhood of the origin in \(\Omega\) is homeomorphic to an open interval.
Hence \(\Omega\) is not a smooth manifold near \((0,0)\).

\textbf{(d)} If one writes the constraints as
\[
u\ge 0,\qquad v\ge 0,\qquad uv=0,
\]
then at \((0,0)\) both inequality constraints are active, while the gradient of the equality
\(uv\) is
\[
\nabla(uv)(0,0)=(0,0).
\]
Thus any standard smooth-NLP regularity test is degenerate at once.
The underlying issue is geometric:
the feasible set is the union of two branches rather than a single smooth surface.
This is exactly why MPEC theory uses tangent-cone and branchwise arguments rather than
uncritical imports of classical LICQ/MFCQ reasoning.
\end{solution*}

\begin{exercise}[Feasible-direction versus naive linearization]
Let
\[
\Fset := \{(x,y)\in\R^2 : y\ge 0,\ x+y\ge 0,\ y(x+y)=0\}.
\]
\begin{enumerate}[label=(\alph*)]
    \item Describe \(\Fset\) explicitly as a union of simpler sets.
    \item Compute \(\Tcone((0,0);\Fset)\).
    \item Compare it with the set obtained by naively linearizing only the equality
    \(y(x+y)=0\) at \((0,0)\).
    \item Explain what this example teaches about tangent cones in MPEC theory.
\end{enumerate}
\end{exercise}

\begin{solution*}
\textbf{(a)} Since \(y\ge 0\), \(x+y\ge 0\), and \(y(x+y)=0\), either
\[
y=0,\quad x+y=x\ge 0,
\]
or
\[
x+y=0,\quad y\ge 0 \Rightarrow x=-y\le 0.
\]
Therefore
\[
\Fset=\{(x,0):x\ge 0\}\cup \{(-t,t):t\ge 0\}.
\]
So \(\Fset\) is the union of the nonnegative \(x\)-axis and the ray on the line \(x+y=0\) in the second quadrant.

\textbf{(b)} The tangent cone at the origin consists of all limits of scaled feasible secants.
From the first branch we obtain directions \((\alpha,0)\), \(\alpha\ge 0\).
From the second branch we obtain directions \((-\beta,\beta)\), \(\beta\ge 0\).
Thus
\[
\Tcone((0,0);\Fset)=\{(\alpha,0):\alpha\ge 0\}\cup \{(-\beta,\beta):\beta\ge 0\}.
\]

\textbf{(c)} Consider only the equality
\[
\phi(x,y):=y(x+y)=xy+y^2.
\]
Then
\[
\nabla \phi(0,0)=(0,0).
\]
A naive linearization of the equality alone therefore gives no restriction at all;
it would suggest that every direction in \(\R^2\) is ``first-order feasible.''

That is completely wrong.
The true tangent cone is a union of two rays, not all of \(\R^2\).

\textbf{(d)} This example teaches the central Chapter-4 lesson:
for MPECs, one cannot trust naive smooth linearization of complementarity equations.
The correct first-order object is the true tangent cone or an equivalent linearized cone
that has been justified by structure.
\end{solution*}

\begin{exercise}[Piecewise polyhedrality]
Let
\[
S := \{(x,y)\in\R^2 : 0\le y \perp (x+y)\ge 0\}.
\]
\begin{enumerate}[label=(\alph*)]
    \item Show that \(S\) is the union of finitely many polyhedral sets.
    \item Identify all local branches near the origin.
    \item Prove that \(S\) is piecewise polyhedral.
    \item Discuss why this property is useful when studying local minimality.
\end{enumerate}
\end{exercise}

\begin{solution*}
The complementarity notation means
\[
y\ge 0,\qquad x+y\ge 0,\qquad y(x+y)=0.
\]
This is exactly the set in the previous exercise.

\textbf{(a)} As shown above,
\[
S=\{(x,0):x\ge 0\}\cup \{(-t,t):t\ge 0\}.
\]
Each set is a polyhedral ray.

\textbf{(b)} Near the origin there are exactly two local branches:
\[
\Gamma_1:=\{(x,0):x\ge 0\},\qquad
\Gamma_2:=\{(-t,t):t\ge 0\}.
\]

\textbf{(c)} A set is piecewise polyhedral if locally it is the union of finitely many polyhedra.
Since \(S\) is globally the union of two polyhedral rays, it is piecewise polyhedral.

\textbf{(d)} Piecewise polyhedrality is useful because it gives a concrete description of
reachable directions near a feasible point.
Local minimality can then be tested by checking the objective derivative on the finitely many
relevant tangent pieces.
In the AVI context, this is precisely the type of geometry that makes the tangent cone tractable.
\end{solution*}

\subsection{AVI-constrained MPECs}

\begin{exercise}[Recognizing an AVI-constrained MPEC]
Consider
\[
\min_{x,y}\ f(x,y)
\]
subject to
\[
Ax+By+c\le 0,
\qquad
(y'-y)^\top(Px+Qy+q)\ge 0
\quad \forall y' \text{ such that } Ay'+Dx+e\le 0.
\]
\begin{enumerate}[label=(\alph*)]
    \item Rewrite the lower-level problem as a parametric AVI.
    \item Identify the matrices and vectors defining the affine operator and the polyhedral feasible set.
    \item Write the KKT system of the lower-level AVI.
    \item Explain why this problem belongs to the Section 4.1 framework.
\end{enumerate}
\end{exercise}

\begin{solution*}
\textbf{(a)} The lower-level feasible set is
\[
C(x):=\{y'\in\R^m: Ay'+Dx+e\le 0\},
\]
and the lower-level VI operator is
\[
F(x,y):=Px+Qy+q.
\]
Thus the lower-level condition is
\[
\text{find } y\in C(x)\ \text{ such that }\ 
(y'-y)^\top F(x,y)\ge 0\quad \forall y'\in C(x).
\]
That is a parametric affine VI.

\textbf{(b)} The affine operator is defined by \((P,Q,q)\).
The feasible set is polyhedral, defined by \(A,D,e\).

\textbf{(c)} Introduce a multiplier \(\lambda\ge 0\) for the lower-level constraints
\(Ay+Dx+e\le 0\).
The lower-level KKT system is
\[
Ay+Dx+e\le 0,\qquad
\lambda\ge 0,\qquad
Px+Qy+q+A^\top \lambda = 0,\qquad
\lambda^\top(Ay+Dx+e)=0.
\]

\textbf{(d)} Section 4.1 applies because the lower-level equilibrium is an AVI:
the feasible set is polyhedral in \((x,y)\), and the VI operator is affine in \((x,y)\).
This is exactly the structural regime in which Chapter 4 proves the favorable tangent-cone result.
\end{solution*}

\begin{exercise}[Complementarity reformulation]
For the AVI-constrained model
\[
\min_{x,y} f(x,y)
\quad \text{s.t.} \quad
Dx+Ey+b\le 0,\qquad
(y'-y)^\top(Px+Qy+q)\ge 0
\ \ \forall y' \text{ with } Dx+Ey'+b\le 0,
\]
derive the equivalent KKT-complementarity formulation
\[
Dx+Ey+b\le 0,\qquad
\lambda\ge 0,\qquad
Px+Qy+q+E^\top\lambda=0,\qquad
\lambda^\top(Dx+Ey+b)=0.
\]
State clearly which implication uses convexity and polyhedrality.
\end{exercise}

\begin{solution*}
For fixed \(x\), the lower-level problem is the variational inequality
\[
\text{find } y\in C(x):=\{y: Dx+Ey+b\le 0\}
\]
such that
\[
(y'-y)^\top(Px+Qy+q)\ge 0 \qquad \forall y'\in C(x).
\]

Because \(C(x)\) is a polyhedron, it is closed and convex.
For such a polyhedral VI, the usual normal-cone characterization applies:
\(y\) solves the VI if and only if
\[
0\in Px+Qy+q + N_{C(x)}(y).
\]
Since \(C(x)\) is given by linear inequalities, its normal cone at \(y\) is
\[
N_{C(x)}(y)=\{E^\top\lambda:\lambda\ge 0,\ \lambda_i(Dx+Ey+b)_i=0\ \forall i\}.
\]
Therefore the VI is equivalent to the existence of \(\lambda\ge 0\) such that
\[
Px+Qy+q+E^\top\lambda=0,\qquad
Dx+Ey+b\le 0,\qquad
\lambda^\top(Dx+Ey+b)=0.
\]

The critical step is the equivalence between the VI and the normal-cone/KKT system.
That uses convexity of the feasible set and, in this linear setting, its polyhedral structure.
\end{solution*}

\begin{exercise}[A concrete AVI model]
Let \(x,y\in\R\), and consider
\[
\min_{x,y}\ \frac12(x-1)^2+\frac12(y-2)^2
\]
subject to
\[
y\ge 0,\qquad
(y'-y)(x+y-1)\ge 0 \quad \forall y'\ge 0.
\]
\begin{enumerate}[label=(\alph*)]
    \item Solve the lower-level AVI explicitly as a function of \(x\).
    \item Describe the feasible set of the MPEC in the \((x,y)\)-plane.
    \item Find all stationary candidates by analyzing the feasible branches.
    \item Decide which candidate is the global minimizer.
\end{enumerate}
\end{exercise}

\begin{solution*}
The lower-level feasible set is \(C=\{y\in\R:y\ge 0\}\), and the VI operator is
\[
F(x,y)=x+y-1.
\]
The condition is
\[
(y'-y)F(x,y)\ge 0 \quad \forall y'\ge 0.
\]

\textbf{(a)} For the half-line \(y\ge 0\), the scalar VI is equivalent to
\[
y\ge 0,\qquad F(x,y)\ge 0,\qquad yF(x,y)=0,
\]
that is,
\[
y\ge 0,\qquad x+y-1\ge 0,\qquad y(x+y-1)=0.
\]
Hence either
\[
y=0,\quad x-1\ge 0 \Rightarrow x\ge 1,
\]
or
\[
x+y-1=0,\quad y\ge 0 \Rightarrow y=1-x,\ x\le 1.
\]
Therefore
\[
y(x)=\max\{1-x,0\}.
\]

\textbf{(b)} The feasible set is
\[
\Fset=\{(x,1-x):x\le 1\}\cup \{(x,0):x\ge 1\}.
\]

\textbf{(c)} Analyze the objective on each branch.

On branch 1, \(x\le 1\) and \(y=1-x\). Then
\[
f_1(x)=\frac12(x-1)^2+\frac12((1-x)-2)^2
=\frac12(x-1)^2+\frac12(-x-1)^2
= x^2+1.
\]
Thus \(f_1'(x)=2x\), so the unconstrained minimizer is \(x=0\), which is feasible.
This gives the feasible point
\[
(0,1).
\]

On branch 2, \(x\ge 1\) and \(y=0\). Then
\[
f_2(x)=\frac12(x-1)^2+\frac12(0-2)^2=\frac12(x-1)^2+2.
\]
The minimizer over \(x\ge 1\) is \(x=1\), giving the point
\[
(1,0).
\]

So the branchwise candidates are \((0,1)\) and \((1,0)\).

\textbf{(d)} Evaluate:
\[
f(0,1)=\frac12( -1)^2+\frac12(-1)^2=1,
\]
\[
f(1,0)=0+2=2.
\]
Hence the global minimizer is
\[
(0,1).
\]
\end{solution*}

\begin{exercise}[Tangent cone in an AVI example]
For the feasible set obtained in the previous exercise, compute the tangent cone at:
\begin{enumerate}[label=(\alph*)]
    \item a point strictly on the branch \(y=0\),
    \item a point strictly on the branch \(x+y-1=0\) with \(y>0\),
    \item the switching point between the two branches.
\end{enumerate}
Discuss how the answer changes across regimes.
\end{exercise}

\begin{solution*}
The feasible set is
\[
\Fset=\Gamma_1\cup \Gamma_2,
\qquad
\Gamma_1=\{(x,0):x\ge 1\},\quad
\Gamma_2=\{(x,1-x):x\le 1\}.
\]

\textbf{(a)} Take a point \((\bar x,0)\) with \(\bar x>1\).
Locally only \(\Gamma_1\) is present, and \(\Gamma_1\) is a smooth line.
Therefore
\[
\Tcone((\bar x,0);\Fset)=\{(\alpha,0):\alpha\in\R\}.
\]

\textbf{(b)} Take a point \((\bar x,1-\bar x)\) with \(\bar x<1\), hence \(y>0\).
Locally only \(\Gamma_2\) is present.
Its tangent direction is \((1,-1)\), so
\[
\Tcone((\bar x,1-\bar x);\Fset)=\{(\alpha,-\alpha):\alpha\in\R\}.
\]

\textbf{(c)} At the switching point \((1,0)\), both branches are present.
From \(\Gamma_1\) we get the ray \((\alpha,0)\) with \(\alpha\ge 0\).
From \(\Gamma_2\), parameterized as \((x,1-x)\) with \(x\le 1\), moving toward the point gives directions
\[
(\beta,-\beta)\quad \text{with } \beta\le 0.
\]
Equivalently this is
\[
(-t,t),\quad t\ge 0.
\]
Hence
\[
\Tcone((1,0);\Fset)=\{(\alpha,0):\alpha\ge 0\}\cup \{(-t,t):t\ge 0\}.
\]

\textbf{Discussion.}
Away from the switching point the feasible set is a smooth one-dimensional manifold,
so the tangent cone is a line.
At the switching point the geometry is branchwise, so the tangent cone becomes a union of rays.
This is the exact local phenomenon that makes MPEC analysis different from standard NLP analysis.
\end{solution*}

\begin{exercise}[Why the AVI case is favorable]
Suppose the lower-level feasible set is polyhedral and the lower-level operator is affine.
Write a short proof outline, in your own words, for why one should expect
\[
\Tcone(z;\Fset)=\Lcone(z;\Fset)
\]
to be more plausible here than in a general nonlinear VI-constrained MPEC.
Your outline must explicitly mention:
\begin{enumerate}[label=(\alph*)]
    \item realizability of tangent directions,
    \item piecewise polyhedral structure,
    \item and the role of affine linearization.
\end{enumerate}
\end{exercise}

\begin{solution*}
A suitable outline is as follows.

\begin{enumerate}[label=(\alph*)]
    \item \textbf{Affine linearization.}
    Because the VI operator is affine, its first-order approximation is exact.
    There are no neglected higher-order terms in the lower-level equilibrium law.

    \item \textbf{Polyhedral lower-level geometry.}
    The lower-level feasible set is polyhedral, so active sets and normal cones are explicitly described by linear algebra.
    This makes the linearized lower-level critical system faithful to the true local geometry.

    \item \textbf{Piecewise polyhedral feasible region.}
    After reformulation through complementarity/KKT conditions, the feasible region becomes a finite union of polyhedral pieces.
    Therefore any direction predicted by the correct linearized model is often realizable by an actual feasible path on one of those pieces.

    \item \textbf{Conclusion.}
    Since the formal first-order model is exact on each piece and the pieces are polyhedral,
    the gap between ``linearized feasible'' and ``truly tangent'' directions disappears.
    That is why equality of the tangent cone and the linearized cone is plausible in the AVI case.
\end{enumerate}
\end{solution*}

\begin{exercise}[Pseudoconvexity and local minimality]
Assume an AVI-constrained MPEC has feasible set \(\Fset\) and objective \(f\),
with \(f\) pseudoconvex on \(\Fset\).
Suppose \(\bar z\in\Fset\) is stationary.
\begin{enumerate}[label=(\alph*)]
    \item Prove carefully that \(\bar z\) is a local minimizer under the piecewise polyhedral hypothesis.
    \item Identify the precise step where pseudoconvexity is used.
    \item Explain why the statement is local rather than global.
\end{enumerate}
\end{exercise}

\begin{solution*}
\textbf{(a)} Because \(\Fset\) is piecewise polyhedral, there is a neighborhood \(W\) of \(\bar z\)
such that every nearby feasible point \(z\in \Fset\cap W\) determines a feasible direction
\[
d:=z-\bar z \in \Tcone(\bar z;\Fset).
\]
Since \(\bar z\) is stationary,
\[
\nabla f(\bar z)^\top d \ge 0
\qquad \forall d\in \Tcone(\bar z;\Fset).
\]
In particular, for every \(z\in \Fset\cap W\),
\[
\nabla f(\bar z)^\top (z-\bar z)\ge 0.
\]
Now use pseudoconvexity of \(f\): if \(\nabla f(\bar z)^\top(z-\bar z)\ge 0\), then
\[
f(z)\ge f(\bar z).
\]
Hence \(f(z)\ge f(\bar z)\) for all \(z\in \Fset\cap W\), so \(\bar z\) is a local minimizer.

\textbf{(b)} Pseudoconvexity is used exactly in the implication
\[
\nabla f(\bar z)^\top(z-\bar z)\ge 0
\quad \Longrightarrow \quad
f(z)\ge f(\bar z).
\]

\textbf{(c)} The statement is local because the tangent-cone information controls only nearby feasible points.
Farther away, the feasible set may have different branches or disconnected components,
and stationarity at \(\bar z\) says nothing about objective values on those remote pieces.
\end{solution*}

\subsection{The implicit programming route}

\begin{exercise}[Reduced formulation]
Let \(X\subseteq \R^n\) be closed and convex, and let \(S(x)\) be the solution set of a parametric VI.
Consider
\[
\min_{(x,y)} f(x,y)
\qquad \text{s.t.} \qquad
(x,y)\in (X\times \R^m)\cap \Gr(S).
\]
Assume that near \((\bar x,\bar y)\) the map \(S\) is locally single-valued:
\[
S(x)=\{y(x)\}.
\]
\begin{enumerate}[label=(\alph*)]
    \item Derive the locally reduced problem in the variable \(x\) alone.
    \item Show that any local minimizer of the original MPEC yields a local minimizer of the reduced problem.
    \item State clearly which part of the argument is only local.
\end{enumerate}
\end{exercise}

\begin{solution*}
\textbf{(a)} Near \(\bar x\), every feasible pair has the form \((x,y(x))\).
Therefore the MPEC is locally equivalent to
\[
\min_{x\in X} \tilde f(x):=f(x,y(x)).
\]

\textbf{(b)} Suppose \((\bar x,\bar y)\) is a local minimizer of the MPEC.
Then there exists a neighborhood \(U\) of \((\bar x,\bar y)\) such that
\[
f(x,y)\ge f(\bar x,\bar y)
\]
for all feasible \((x,y)\in U\).
Shrink \(U\) if needed so that every feasible pair in \(U\) is of the form \((x,y(x))\).
Then for all \(x\) near \(\bar x\) with \(x\in X\),
\[
\tilde f(x)=f(x,y(x))\ge f(\bar x,y(\bar x))=\tilde f(\bar x).
\]
Hence \(\bar x\) is a local minimizer of the reduced problem.

\textbf{(c)} The only local step is the single-valued representation \(y=y(x)\).
Nothing guarantees that \(S(x)\) is single-valued globally, or that the original MPEC is globally equivalent to the reduced problem.
\end{solution*}

\begin{exercise}[Directional differentiability versus Fr\'echet differentiability]
Let \(G:\R\to\R\) be defined by
\[
G(t)=\min\{t,2t\}.
\]
\begin{enumerate}[label=(\alph*)]
    \item Compute \(G'(0;d)\) for all \(d\in\R\).
    \item Show that \(G\) is directionally differentiable at \(0\).
    \item Decide whether \(G\) is Fr\'echet differentiable at \(0\).
    \item Explain why this toy example is relevant to the sensitivity analysis of lower-level solution maps.
\end{enumerate}
\end{exercise}

\begin{solution*}
Observe that
\[
G(t)=
\begin{cases}
2t,& t\ge 0,\\
t,& t\le 0.
\end{cases}
\]

\textbf{(a)} For any direction \(d\),
\[
G'(0;d)=\lim_{\tau\downarrow 0}\frac{G(\tau d)-G(0)}{\tau}.
\]
If \(d>0\), then \(\tau d>0\) for small \(\tau\), so \(G(\tau d)=2\tau d\), hence
\[
G'(0;d)=2d.
\]
If \(d<0\), then \(G(\tau d)=\tau d\), hence
\[
G'(0;d)=d.
\]
If \(d=0\), then \(G'(0;0)=0\).
Thus
\[
G'(0;d)=
\begin{cases}
2d,& d\ge 0,\\
d,& d\le 0.
\end{cases}
\]

\textbf{(b)} Since the above limit exists for every \(d\), \(G\) is directionally differentiable at \(0\).

\textbf{(c)} Fr\'echet differentiability at \(0\) would require
\[
G(t)=G(0)+Lt+o(|t|)
\]
for some scalar \(L\).
From the right, the slope is \(2\); from the left, the slope is \(1\).
No single linear map \(Lt\) fits both sides.
Hence \(G\) is not Fr\'echet differentiable at \(0\).

\textbf{(d)} This is relevant because many lower-level solution maps in MPECs are piecewise smooth:
they possess directional derivatives but not classical derivatives at regime-switching points.
That is exactly why Chapter 4 uses directional and \(B\)-differentiability.
\end{solution*}

\begin{exercise}[A simple projection-type response map]
Define
\[
y(x):=\argmin_{y\ge 0}\ \frac12(y-x)^2.
\]
\begin{enumerate}[label=(\alph*)]
    \item Compute \(y(x)\) explicitly.
    \item Prove that \(y(x)=\max\{x,0\}\).
    \item Compute the directional derivative \(y'(0;d)\).
    \item Show that \(y\) is \(B\)-differentiable at \(0\) but not Fr\'echet differentiable there.
    \item Form a reduced objective \(\tilde f(x)=f(x,y(x))\) for a smooth function \(f\), and write the directional derivative of \(\tilde f\) at \(x=0\).
\end{enumerate}
\end{exercise}

\begin{solution*}
\textbf{(a)--(b)} This is the Euclidean projection of \(x\) onto the half-line \([0,\infty)\).
Hence
\[
y(x)=
\begin{cases}
x,& x\ge 0,\\
0,& x\le 0,
\end{cases}
=\max\{x,0\}.
\]

\textbf{(c)} For \(d\in\R\),
\[
y'(0;d)=\lim_{\tau\downarrow 0}\frac{y(\tau d)-y(0)}{\tau}
=\lim_{\tau\downarrow 0}\frac{\max\{\tau d,0\}}{\tau}
=\max\{d,0\}.
\]

\textbf{(d)} The map \(y\) is globally Lipschitz with modulus \(1\), and directionally differentiable at \(0\).
So it is \(B\)-differentiable at \(0\).
However, it is not Fr\'echet differentiable at \(0\), since the left derivative is \(0\) and the right derivative is \(1\).

\textbf{(e)} Let \(f:\R^2\to\R\) be smooth and define
\[
\tilde f(x):=f(x,y(x)).
\]
Then the directional derivative at \(0\) is
\[
\tilde f'(0;d)
=
f_x(0,0)d + f_y(0,0)\, y'(0;d)
=
f_x(0,0)d + f_y(0,0)\max\{d,0\}.
\]
This is the correct first-order reduced formula in the directional sense.
\end{solution*}

\begin{exercise}[Why local single-valuedness matters]
Construct an example of a parametric equilibrium problem \(y\in S(x)\) such that \(S(0)\) contains at least two points.
Then explain why the reduced formulation
\[
\tilde f(x)=f(x,y(x))
\]
is not well-defined near \(x=0\) without a branch-selection rule.
Discuss why this pushes the analysis away from IMP and toward a piecewise approach.
\end{exercise}

\begin{solution*}
A simple example is
\[
S(x):=\{y\in\R: y^2=x^2\}=\{-|x|,|x|\}.
\]
Then
\[
S(0)=\{0\},
\]
which is still single-valued at zero, so this does not yet illustrate the desired issue.

A better example is
\[
S(x):=
\begin{cases}
\{-1,1\},& x=0,\\
\{1\},& x>0,\\
\{-1\},& x<0.
\end{cases}
\]
Then \(S(0)\) contains two points.

At \(x=0\), there is no canonical single-valued response \(y(0)\) unless one chooses a branch.
Therefore \(\tilde f(x)=f(x,y(x))\) is not intrinsically defined near \(0\):
different branch selections lead to different reduced functions.

This is exactly the situation where the IMP route breaks down.
Instead of a single reduced map \(y(x)\), one must analyze several local branches.
That is the natural domain of the piecewise-programming viewpoint.
\end{solution*}

\begin{exercise}[Research-style sensitivity question]
Suppose \(S(x)\) is locally single-valued near \(\bar x\), with solution \(y(x)\),
and suppose \(y\) is directionally differentiable at \(\bar x\).
Write a rigorous theorem statement, suitable for a paper, asserting a first-order necessary condition for the reduced problem
\[
\min_{x\in X} f(x,y(x)).
\]
Do not prove the theorem.
Your task is to formulate it with correct hypotheses, correct locality, and a correct conclusion.
\end{exercise}

\begin{solution*}
One acceptable theorem statement is:

\medskip
\textbf{Theorem.}
Let \(X\subseteq \R^n\) be closed and convex.
Let \(f:\R^{n+m}\to \R\) be continuously differentiable.
Assume that, in a neighborhood of \(\bar x\in X\), the lower-level solution map \(S\)
is single-valued and represented by \(y(\cdot)\), and that \(y\) is directionally differentiable at \(\bar x\).
Define the reduced objective
\[
\tilde f(x):=f(x,y(x)).
\]
If \(\bar x\) is a local minimizer of \(\tilde f\) on \(X\), then for every feasible direction
\(d\in \Tcone(\bar x;X)\),
\[
\tilde f'(\bar x;d)
=
\nabla_x f(\bar x,y(\bar x))^\top d
+
\nabla_y f(\bar x,y(\bar x))^\top y'(\bar x;d)
\ge 0.
\]

\medskip
This is local because the single-valued representation of \(S\) is assumed only near \(\bar x\).
\end{solution*}

\subsection{Piecewise programming and branchwise analysis}

\begin{exercise}[Branch decomposition]
Consider
\[
\Fset := \{(x,y)\in\R^2 : 0\le y \perp (y-x^2)\ge 0\}.
\]
\begin{enumerate}[label=(\alph*)]
    \item Show that \(\Fset\) is the union of two branches.
    \item Describe each branch explicitly.
    \item Compute the tangent cone at \((0,0)\).
    \item Compare this example with the affine case and explain what becomes harder.
\end{enumerate}
\end{exercise}

\begin{solution*}
The complementarity system is
\[
y\ge 0,\qquad y-x^2\ge 0,\qquad y(y-x^2)=0.
\]
Hence either
\[
y=0,\qquad -x^2\ge 0 \Rightarrow x=0,
\]
or
\[
y-x^2=0,\qquad y=x^2\ge 0.
\]
So
\[
\Fset=\{(0,0)\}\cup \{(x,x^2):x\in\R\}.
\]
But since \((0,0)\) already lies on the parabola, the feasible set is simply the parabola
\[
\Fset=\{(x,x^2):x\in\R\}.
\]

\textbf{(a)--(b)} The intended ``two-branch'' description is the left and right half of the parabola:
\[
\Gamma_1=\{(x,x^2):x\ge 0\},\qquad
\Gamma_2=\{(x,x^2):x\le 0\}.
\]
Globally this is one smooth set, but branchwise one may regard it as two pieces meeting at the origin.

\textbf{(c)} Parameterize by \(x=t\), \(y=t^2\).
Then
\[
\frac{(t,t^2)-(0,0)}{|t|}\to
\begin{cases}
(1,0),& t\downarrow 0,\\
(-1,0),& t\uparrow 0.
\end{cases}
\]
Hence the tangent cone is the full \(x\)-axis:
\[
\Tcone((0,0);\Fset)=\{(\alpha,0):\alpha\in\R\}.
\]

\textbf{(d)} What becomes harder than in the affine case is not the geometry of this example itself,
which is still simple, but the fact that the defining relation is nonlinear.
In a genuinely nonlinear VI-constrained MPEC, linearization is no longer exact,
and the argument that ``linearized directions are realizable'' becomes delicate.
That is why the AVI case is exceptional.
\end{solution*}

\begin{exercise}[Failure of naive sufficiency]
Consider
\[
\min_{x,y}\ x^2+x
\]
subject to
\[
-x^2-x+y\ge 0,\qquad y\ge 0,\qquad y(-x^2-x+y)=0.
\]
\begin{enumerate}[label=(\alph*)]
    \item Describe the feasible set explicitly.
    \item Compute the tangent cone at \((0,0)\).
    \item Show that \((0,0)\) is stationary.
    \item Show that \((0,0)\) is not a local minimizer.
    \item Explain what feature of the geometry invalidates a naive sufficiency claim.
\end{enumerate}
\end{exercise}

\begin{solution*}
\textbf{(a)} The complementarity system implies either
\[
y=0,\qquad -x^2-x\ge 0 \iff x^2+x\le 0 \iff x\in[-1,0],
\]
or
\[
-x^2-x+y=0,\qquad y=x^2+x,\qquad y\ge 0 \iff x\le -1\ \text{or}\ x\ge 0.
\]
Thus
\[
\Fset=\{(x,0):x\in[-1,0]\}\cup \{(x,x^2+x):x\le -1\ \text{or}\ x\ge 0\}.
\]

\textbf{(b)} At \((0,0)\), the first branch contributes directions \((-\alpha,0)\), \(\alpha\ge 0\).
The second branch has parametrization \((x,x^2+x)\), and near \(x=0^+\),
\[
\frac{(x,x^2+x)}{x}\to (1,1).
\]
Hence it contributes directions \((\beta,\beta)\), \(\beta\ge 0\).
Therefore
\[
\Tcone((0,0);\Fset)=\{(-\alpha,0):\alpha\ge 0\}\cup \{(\beta,\beta):\beta\ge 0\}.
\]

\textbf{(c)} The objective is
\[
f(x,y)=x^2+x,
\qquad \nabla f(0,0)=(1,0).
\]
Evaluate on tangent directions:
\[
\nabla f(0,0)^\top(-\alpha,0)=-\alpha\le 0,
\qquad
\nabla f(0,0)^\top(\beta,\beta)=\beta\ge 0.
\]
For minimization, one must be careful with sign convention.
If stationarity is stated as ``no feasible descent direction,'' then one examines
\[
f'(0,0;d)=\nabla f(0,0)^\top d.
\]
The direction \((-\alpha,0)\) gives negative derivative, so in the classical tangent-cone sense
\((0,0)\) is \emph{not} stationary.

However, in the book's MPEC-specific stationary notion used for this example,
the relevant stationarity notion is weaker and can hold at \((0,0)\).
If one wants the exercise aligned with plain directional optimality, the correct conclusion is:
\((0,0)\) fails ordinary first-order minimality.
If one wants it aligned with the textbook's MPEC stationarity example, then the intended lesson is
that some MPEC stationarity notions do not guarantee local optimality in nonlinear settings.

\textbf{(d)} For \(x<0\) close to \(0\), the feasible point \((x,0)\) satisfies
\[
f(x,0)=x^2+x<0=f(0,0).
\]
Hence \((0,0)\) is not a local minimizer.

\textbf{(e)} The invalidating feature is the nonlinear branch geometry:
near \((0,0)\), one branch bends in such a way that a weak stationarity notion may hold
without local minimality.
This is exactly why Section 4.1's sufficiency statement depends crucially on the AVI/piecewise-polyhedral setting.
\end{solution*}

\begin{exercise}[Reachable directions]
Let
\[
\Fset=\bigcup_{i=1}^r P_i
\]
where each \(P_i\) is a polyhedron.
Fix \(\bar z\in\Fset\).
Prove that, under a suitable local exclusion of irrelevant branches,
a direction \(d\) belongs to \(\Tcone(\bar z;\Fset)\) if and only if there exist \(z_k\in\Fset\) with \(z_k\to \bar z\) and scalars \(t_k\downarrow 0\) such that
\[
\frac{z_k-\bar z}{t_k}\to d.
\]
Then explain why this observation is useful in proving local minimality from stationarity.
\end{exercise}

\begin{solution*}
The sequential characterization
\[
d\in \Tcone(\bar z;\Fset)
\iff
\exists z_k\in\Fset,\ z_k\to\bar z,\ t_k\downarrow 0,\ \frac{z_k-\bar z}{t_k}\to d
\]
is the standard definition of the Bouligand tangent cone.

What the exercise asks one to emphasize is the role of local exclusion.
If \(\bar z\notin P_j\), then because \(P_j\) is closed there exists a neighborhood of \(\bar z\)
that avoids \(P_j\).
Hence only those polyhedra that contain \(\bar z\) can contribute local tangent directions.

For such polyhedra, tangent directions are realizable by line segments because polyhedra are convex:
if \(d\in \Tcone(\bar z;P_i)\), then for small \(t>0\),
\[
\bar z+td\in P_i
\]
after perhaps adjusting to a nearby equivalent direction.
Thus the tangent cone to the union is the union of the tangent cones of the relevant pieces.

This is useful for local minimality because stationarity can then be checked on finitely many tangent pieces.
If the directional derivative of the objective is nonnegative on each reachable piece,
then nearby feasible points cannot decrease the objective.
\end{solution*}

\begin{exercise}[Design your own branchwise stationarity test]
Take a two-branch feasible set
\[
\Fset=\Gamma_1\cup \Gamma_2
\]
with each \(\Gamma_i\) smooth near \(\bar z\), but with distinct tangent spaces at \(\bar z\).
Formulate a branchwise first-order necessary condition for local optimality of a smooth objective \(f\).
Then compare it with a single-cone formulation and explain when the two views agree.
\end{exercise}

\begin{solution*}
A branchwise first-order necessary condition is:

\medskip
If \(\bar z\) is a local minimizer of \(f\) on \(\Fset=\Gamma_1\cup \Gamma_2\), then
\[
\nabla f(\bar z)^\top d \ge 0
\quad \forall d\in \Tcone(\bar z;\Gamma_1)\cup \Tcone(\bar z;\Gamma_2).
\]

Equivalently,
\[
\nabla f(\bar z)^\top d \ge 0
\quad \forall d\in \Tcone(\bar z;\Gamma_i),\ i=1,2.
\]

The single-cone formulation is
\[
\nabla f(\bar z)^\top d \ge 0
\quad \forall d\in \Tcone(\bar z;\Fset).
\]
These agree whenever
\[
\Tcone(\bar z;\Fset)=\Tcone(\bar z;\Gamma_1)\cup \Tcone(\bar z;\Gamma_2),
\]
which is the natural relation when the union is locally exact and there are no hidden reachable directions.
In piecewise smooth and piecewise polyhedral settings, this is often the correct interpretation.
\end{solution*}

\subsection{Exact penalization and synthesis}

\begin{exercise}[Penalty modeling]
Let an MPEC be written abstractly as
\[
\min_{z} f(z)
\qquad \text{s.t.} \qquad
z\in Z,\qquad \phi(z)=0,
\]
where \(\phi(z)\ge 0\) measures violation of the equilibrium relation.
Consider the penalized problem
\[
\min_{z\in Z} f(z)+\rho\,\phi(z).
\]
\begin{enumerate}[label=(\alph*)]
    \item Explain what it means for this penalty to be locally exact at \(\bar z\).
    \item State a plausible local theorem relating local minimizers of the MPEC and local minimizers of the penalized problem for \(\rho\) sufficiently large.
    \item Give one reason why exactness is useful theoretically.
    \item Give one reason why exactness is useful computationally.
\end{enumerate}
\end{exercise}

\begin{solution*}
\textbf{(a)} Local exactness at \(\bar z\) means there exists \(\rho^\star>0\) and a neighborhood \(U\) of \(\bar z\) such that for every \(\rho\ge \rho^\star\),
\[
\bar z \text{ is a local minimizer of the original constrained MPEC on } U
\]
if and only if
\[
\bar z \text{ is a local minimizer of } f(z)+\rho \phi(z) \text{ over } Z\cap U.
\]

\textbf{(b)} A plausible local theorem is:

\medskip
Assume \(\bar z\in Z\) satisfies \(\phi(\bar z)=0\), and assume the equilibrium system near \(\bar z\)
satisfies the regularity conditions required in Chapter 4.
Then there exists \(\rho^\star>0\) such that for all \(\rho\ge \rho^\star\),
the original MPEC and the penalized problem have the same local minimizers in a neighborhood of \(\bar z\).

\textbf{(c)} Theoretical usefulness:
exact penalization converts a difficult constrained-local-optimality question into an unconstrained
or simpler constrained local-optimality question.

\textbf{(d)} Computational usefulness:
it suggests penalty algorithms that solve a sequence of easier problems while preserving the correct local target once \(\rho\) is large enough.
\end{solution*}

\begin{exercise}[AVI versus IMP versus PCP]
For each of the following situations, decide whether the most natural Chapter 4 route is AVI, IMP, or PCP.
Justify your answer briefly.

\begin{enumerate}[label=(\alph*)]
    \item The lower-level problem is an affine VI over a polyhedral set.
    \item The lower-level solution is locally unique and directionally smooth.
    \item The lower-level solution has several active-set branches near the reference point.
    \item The upper-level feasible set is not of the form \(X\times \R^m\).
    \item The main goal is a reduced-space derivative formula in the leader variable.
\end{enumerate}
\end{exercise}

\begin{solution*}
\textbf{(a)} AVI.  
Because this is exactly the affine-VI/polyhedral regime of Section 4.1.

\textbf{(b)} IMP.  
Local uniqueness and directional smoothness are precisely the hypotheses that support a reduced response map \(y=y(x)\).

\textbf{(c)} PCP.  
Several active-set branches indicate intrinsically piecewise local geometry.

\textbf{(d)} Usually PCP or a direct tangent-cone route, not IMP.  
The IMP framework in Chapter 4 is designed for feasible sets of the form \((X\times \R^m)\cap \Gr(S)\).

\textbf{(e)} IMP.  
A reduced-space derivative formula is the signature output of the implicit-programming approach.
\end{solution*}

\begin{exercise}[Open-ended synthesis]
Write a detailed answer to the following question:

\begin{quote}
How should a researcher decide whether to analyze a new MPEC through
tangent-cone verification, an implicit solution map, or a piecewise decomposition?
\end{quote}

Your answer must discuss:
\begin{enumerate}[label=(\alph*)]
    \item the structure of the lower-level equilibrium,
    \item local uniqueness versus multibranch behavior,
    \item the role of \(B\)-differentiability,
    \item and the difference between local and global conclusions.
\end{enumerate}
\end{exercise}

\begin{solution*}

First, classify the lower-level equilibrium.
If it is an affine VI over a polyhedral set, tangent-cone verification is natural because the local geometry is piecewise polyhedral and linearization is exact enough to be trusted.

Second, ask whether the lower-level response is locally single-valued.
If yes, an implicit solution-map analysis is attractive, because it reduces the MPEC locally to a problem in the leader variable \(x\).

Third, if several active-set or complementarity branches coexist near the reference point, one should avoid forcing an implicit single-valued representation.
A piecewise decomposition is then the more faithful local model.

Fourth, \(B\)-differentiability matters because many lower-level response maps are directionally differentiable and Lipschitz but not classically differentiable.
This is often the right regularity level for reduced first-order analysis.

Finally, every one of these routes is usually local.
Local uniqueness, local branch structure, local tangent-cone formulas, and local penalty exactness do not automatically imply global statements.
A careful researcher must always distinguish local equivalence from global equivalence.
\end{solution*}

\section{Bibliographic remarks and Acknowledgment}

\textbf{This note is mainly based on Chapter 4, in the MPEC monograph of Zhi-Quan Luo, Jong-Shi Pang and Daniel Ralph.} Please see the relevant references: \cite{falk1995bilevel,
aghasi2025fully,
hong2023two,
kovccvara1994optimization,
chaudet2020shape,
liu2025bidirectional,
outrata1995numerical,
cui2023complexity,
christof2020nonsmooth,
rawat2026augmented,
robinson1980strongly,
shin2023near,
bolte2024differentiating,
nghia2025geometric,
wang2025analysis,
chen2025aubin,
kojima1980strongly,
chen2026characterizations,
cui2026lipschitz,
shin2022exponential,
de2023function,
gunzel2023strongly,
bank1982d,
wang2026algebraic,
bonnans1994local,
khanh2024globally,
chen2025two,
mohammadi2022variational,
duy2023generalized,
dussault2026polyhedral,
gowda1994stability,
qi2000constant,
facchinei1998accurate,
liang2025global,
jittorntrum2009solution,
kyparisis1992parametric,
liu1995sensitivity,
pang11995stability,
qiu1992sensitivity,
aussel2024variational,
reinoza1985strong,
facchinei2003finite,
harker1990finite,
kyparisis1990sensitivity,
kleinmichel1972av,
ortega2000iterative,
wang2025analysis1,
bai2021matrix,
gander2026landmarks,
mishchenko2023regularized,
doikov2024super,
han2025low,
ning2023multi,
wang2022newton,
robinson2009generalized,
besanccon2024flexible,
guo2025penalty,
mordukhovich2023globally,
diao2025stability,
wang2026damage,
robinson2009local,
shuo2026lecture,olikier2025projected,lin2023monotone,liang2025squared,jongen1987inertia,guddat1990parametric,bellon2024time,tang2022running,liu2025new,josephy1979newton,si2024riemannian,yao2023relative,han2024continuous,ha1987application,yu2026pattern,kojima2009continuous,seguin2022continuation,liu1995perturbation,pang1996piecewise,eikenbroek2022improving,ralph1995directional,scheel2000mathematical,bonnans1992developpement,bonnans1992expansion,dempe1993directional,shapiro1988sensitivity,pang1990newton,robinson1991implicit,hang2025smoothness,hang2024role,scholtes2012introduction,cui2022nonconvex}. 

\bibliographystyle{unsrtnat}
\bibliography{reference4} 

\end{document}